\providecommand{\href}[2]{#2}

\documentclass[12pt]{amsart}

\usepackage{graphicx}
\usepackage{amsthm}
\usepackage{amsmath, amssymb}
\swapnumbers

\theoremstyle{plain}
\newtheorem{Thm}{Theorem}

\newtheorem{Coro}[Thm]{Corollary}

\newtheorem{Claim}[Thm]{Claim}

\theoremstyle{definition}
\newtheorem{Def}[Thm]{Definition}

\begin{document}
\begin{abstract}  Let $S$ be a triangulated $2$-sphere with fixed triangulation $T$.  We apply the methods of thin position from knot theory to obtain a simple version of the three geodesics theorem for the 2-sphere  \cite{L-S}.  In general these three geodesics may be unstable, corresponding, for example,  to the three equators of an ellipsoid.   Using a piece-wise linear approach, we show that we can usually find at least three stable geodesics.    

\end{abstract}
\title{ Finding geodesics in a triangulated $2$-sphere}

\maketitle

 \author Abigail Thompson \footnote{Supported in part by the National Science Foundation.}

\section{Introduction}

Let $S$ be a triangulated $2$-sphere with fixed triangulation $T$.  We assume $T$ is a simplicial complex, see \cite{H}.  In particular we may assume that the 1-skeleton of $T$ contains no loops or multiple edges.   We apply the methods of thin position from knot theory to obtain a simple version of the three geodesics theorem for the 2-sphere  \cite{L-S}.  Using this piece-wise linear approach we can go further, and strengthen the result to find at least three stable (PL) geodesics, unless the triangulation is the tetrahedral triangulation or the ``double tetrahedral" triangulation.
 
\subsection{Outline of the paper}
In section 2, we define stable and unstable geodesics, and thin position for a triangulation of the 2-sphere.  We prove the basic result that a thin triangulation naturally yields geodesics corresponding to stable and unstable geodesics.    In section 3 we define bridge position for a triangulation, analogous to bridge position for a knot in the 3-sphere.   We use a result of H. Whitney on the existence of Hamiltonian cycles to examine the relation between thin position and bridge position for a triangulation, and conclude that thin position is the same as bridge position only in the case of the tetrahedral triangulation.   In sections 4 and 5 we pursue this idea, and use it to obtain a relatively simple version of the three geodesics theorem, in which the three geodesics are allowed to be either stable or unstable.   Finally in section 6 we refine our analysis of a thin triangulation to obtain, with two exceptions, the existence of three stable geodesics for a triangulated 2-sphere.

 \section{Width of a triangulation}

Let $S$ be a triangulated $2$-sphere with fixed triangulation $T$.

\begin{Def}

Let $P={e_1,e_2,....,e_k}$ be an imbedded cycle in the edges of $T$.   Let $T_j$ be  a triangle in $T$ that intersects $P$ in exactly one or exactly two (necessarily adjacent) edges in $P$.   A {\it local move on $P$} replaces the one or two edges of $T_j$ with two or one edges of $T_j$, yielding another imbedded cycle $Q$ with either one more or one fewer edges than $P$.    Call the first kind of move a {\it shortening} of $P$; the second a {\it lengthening}.

\end{Def}

\begin{Def}

Let $P={e_1,e_2,....,e_k}$ be an imbedded cycle in the edges of $T$.    $P$ is a {\it stable geodesic}  if it allows no shortening moves and $P$ is not the boundary of a triangle.

\end{Def}

\begin{figure}[h]
    \centering
    \includegraphics[width=0.8\textwidth]{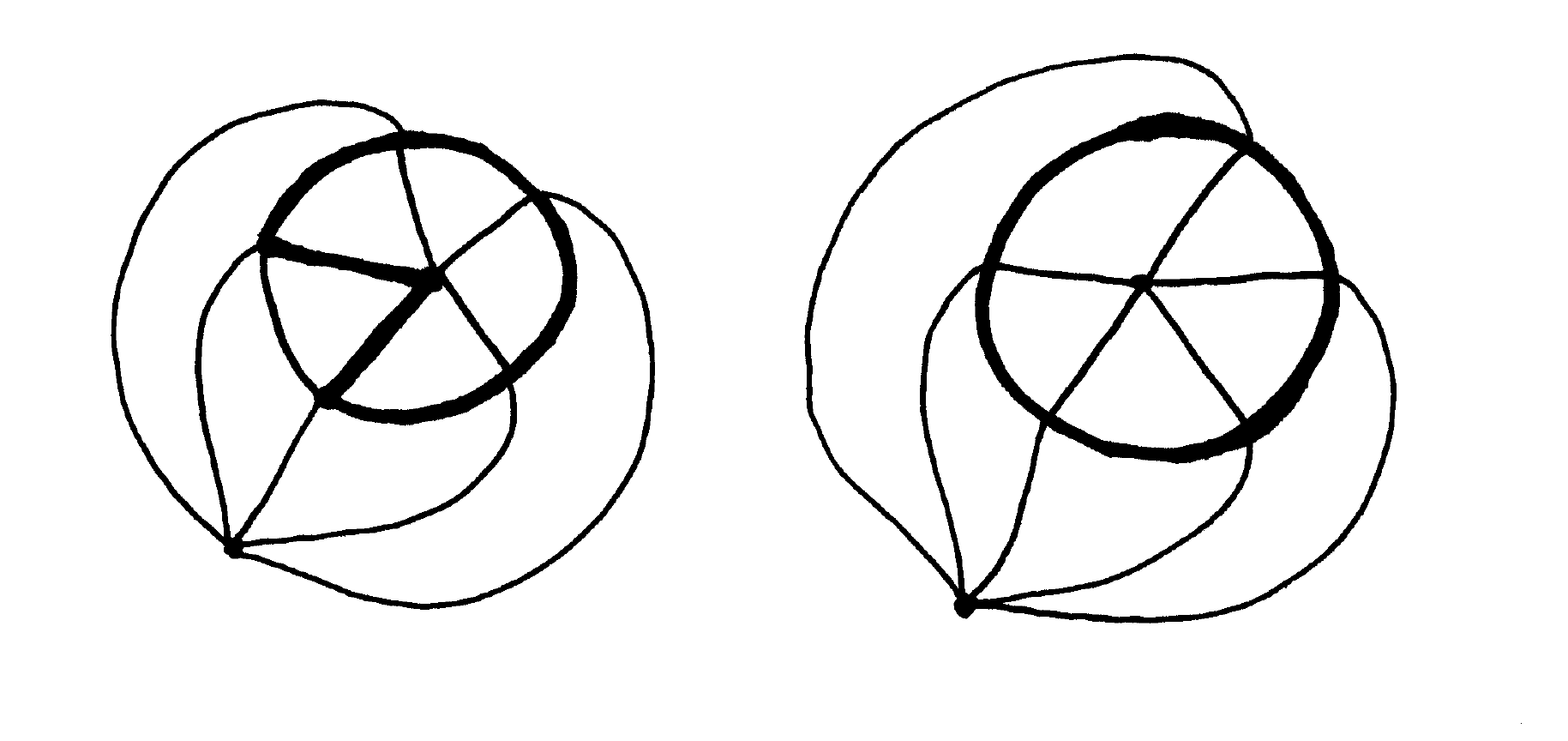}
    \caption{unstable vs. stable geodesic}
    \label{fig1}
\end{figure}

\begin{Def}

Let $P={e_1,e_2,....,e_k}$ be an imbedded cycle in the edges of $T$.   $P$ divides $S$ into two disks, $D_1$ and $D_2$.   Suppose $P$ has two shortening moves, one in $D_1$ across $T_1$  and one in $D_2$ across $T_2$. Suppose further that for every such pair, $T_1$ and $T_2$ intersect in an edge $e_3$ contained in $P$.   Notice that this intersection prevents $P$ from shortening to both sides simultaneously.   We call such a $P$  an {\it unstable geodesic} (see Figure \ref{fig1}).
\end{Def}

\begin{Def}

Since triangulations of the $2$-sphere are shellable \cite{B}, we can choose an order $O(T)$ for the triangles of $T$, $T_1, T_2, .....,T_n$ so that 
$$
I_k=T_1\cup T_2\cup......\cup T_k
$$
is homeomorphic to a disk for $k<n$.    Call such an order {\it good}.  We assume for the remainder of this paper that a specified order for a given $T$ is good.

\end{Def}

\begin{Def}

Let $O(T)$ be an ordering of $T$.   Call the number of vertices of $T$ in the boundary of $I_k$ the {\it length of $\partial(I_k$)}, and denote it $|\partial(I_k)|$.   Notice that in a good ordering, the addition of each successive triangle either increases the length of the boundary of the disk by exactly one or reduces it by exactly one.    A {\it local maximum} of the ordered list $(|\partial I_1|, |\partial I_2|,....,|\partial I_{n-1}|)$ is a value $|\partial I_j|$ such that 
$$
|\partial I_{j-1}|<|\partial I_j|>|\partial I_{j+1}|, 
$$
$j=1,...,n-1.$
 A {\it local minimum} of the ordered list $(|\partial I_1|, |\partial I_2|,....,|\partial I_{n-1}|)$ is a value $|\partial I_j|$ such that 
 $$
 |\partial I_{j-1}|>|\partial I_j|<|\partial I_{j+1}|,
 $$
 $ j=1,...,n-1$.  We say $T$ with order $O(T)$ is in {\it bridge position} if the ordered list $(|\partial I_1|, |\partial I_2|,....,|\partial I_{n-1}|)$ has a single local maximum and no local minima (see Figure \ref{fig2}). 
 
 \begin{figure}[h]
    \centering
    \includegraphics[width=0.8\textwidth]{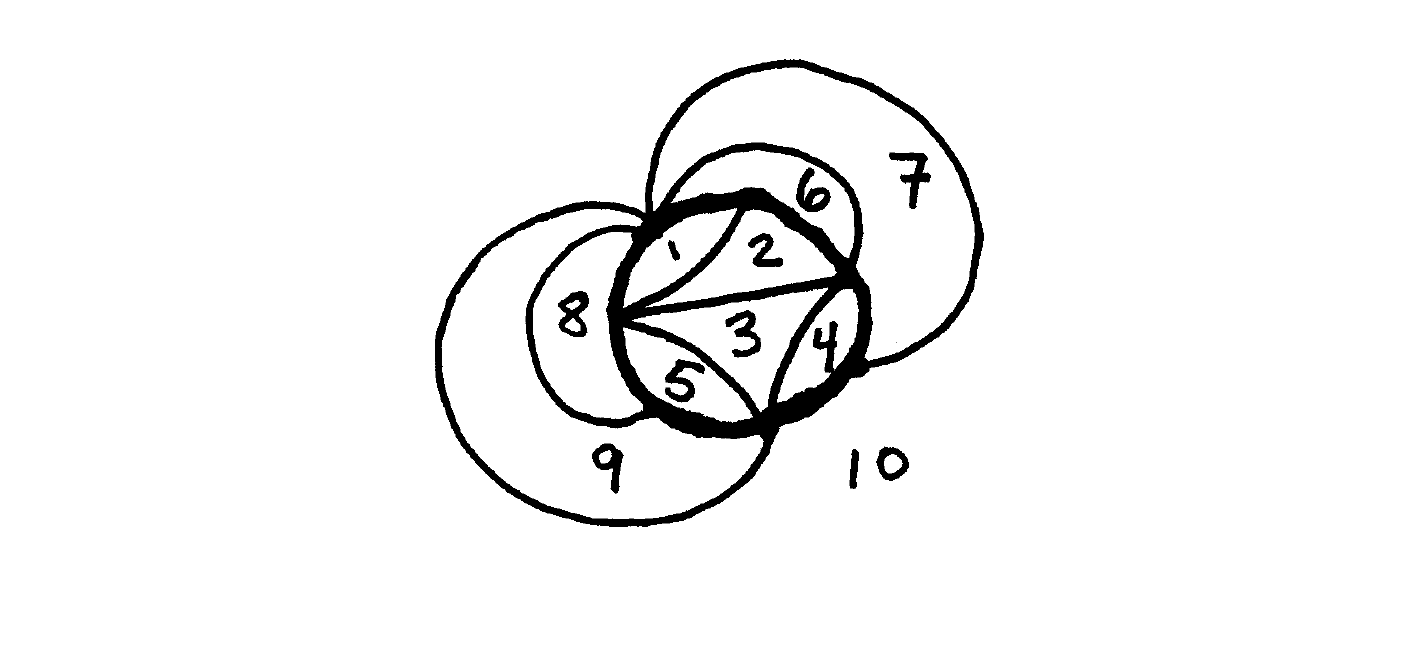}
    \caption{bridge position; local max at $\partial I_5$} 
    \label{fig2}
\end{figure}

 The {\it width of $O(T)$}, {$w_O(T)$}, is the list of local maxima of  \\
$(|\partial I_1|, |\partial I_2|,....,|\partial I_{n-1}|)$, lexicographically ordered.     The {\it width of $(T)$}, $w(T)$,  is the minimum over all such lists, lexicographically ordered.

 We say that $T$ with order $O(T)$ is in {\it thin position} if $O(T)$ realizes the width of $T$.   A local maximum (minimum) {\it corresponds to} the cycle which is the boundary of the disk $\partial I_j$.   

\end{Def} 

\begin{Thm} If $T$ with order $O(T)$ is in thin position, then 
\begin{enumerate}
\item the cycles corresponding to the local maxima are unstable geodesics 
\item the cycles corresponding to the local minima are stable geodesics.
\end{enumerate}
\end{Thm}

\medskip

{\bf Proof}

We start with two technical claims.   We introduce the dual graph of $T$, $\Gamma_T$, which is useful when analyzing how $|\partial(I_k)|$ changes when a triangle is added or removed.    

\begin{Claim} 

Let $D$ be a triangulated disk with a (good) ordering $T_1, T_2,...,T_m$ such that $(|\partial I_1|<|\partial I_2|<....<|\partial I_{m}|)$.  Let $\Gamma_T$ be the dual graph in $D$.   Then $\Gamma_T$ is a tree.  
\end{Claim}

{\bf Proof}

We can build $\Gamma_T$  following the order on $T$.   Since $|\partial(I_k)|$ is strictly increasing as $k$ increases from $1$ to $m$,  as $\Gamma_T$ is built each new vertex must have degree one, hence $\Gamma_T$ is a tree.  

\begin{Claim} 

Let $D$ be a triangulated disk with a (good) ordering $O(T)$, $T_1, T_2,...,T_m$.  Suppose there is a shortening move for $\partial{D}$ across $T_i$.    Then the ordering  $O_*(T)$ given by $T_1,..,T_{i-1},T_{i+1},..,T_m,T_i$ is also a good ordering, and $w_{O_*}(T)\leq{w_O(T)}$.
\end{Claim}

{\bf Proof}

\begin{figure}[h]
    \centering
    \includegraphics[width=0.8\textwidth]{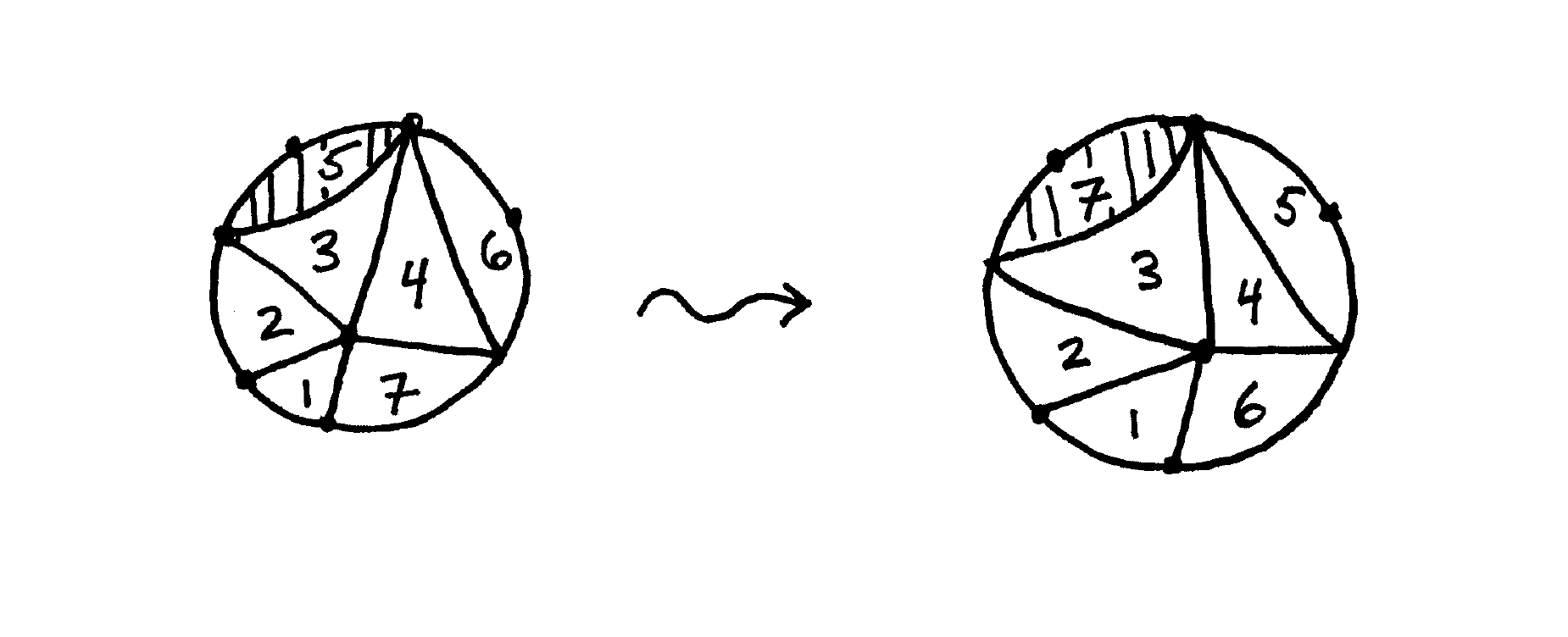}
    \caption{re-ordering shaded triangle}
    \label{fig3}
\end{figure}

Note that since there is a shortening move for $\partial{D}$ across $T_i$,  $T_i$ corresponds to a valence one vertex in $\Gamma_T$.   Thus the homeomorphism type of $T_1\cup{T_2}\cup..\cup {T_{i-1}}\cup{T_{i+1}}\cup..\cup{T_j}$ is the same as that of $I_j$, hence the ordering $T_1,..,T_{i-1},T_{i+1},..,T_m,T_i$ is still good.     Each local maximum in the list $(|\partial I_1|, |\partial I_2|,....,|\partial I_{m}|)$ is either unchanged or reduced by one when the addition of $T_i$ is delayed to the last step.    Call this the {\it re-ordering principle} (see Figure \ref{fig3}).   

Note that the re-ordering principle applies more generally.   Suppose $T_1, T_2,...,T_m$ is a triangulation of a planar region $P$ and  suppose there is a shortening move across $T_i$ for a boundary curve $C$ of $P$.   We can define the width of this ordering for $P$ as before.    By the argument above, each local maximum in the list $(|\partial I_1|, |\partial I_2|,....,|\partial I_{m}|)$ is either unchanged or reduced by one when the addition of $T_i$ is delayed, and the homeomorphism type of the region at each stage is unchanged if the addition of $T_i$ is delayed, so that we may assume $i=m$.     

\bigskip
Suppose $T$ with order $O(T)$ is in thin position, and suppose $C$ is a  cycle corresponding to a local maximum $|\partial I_{j}|$.   

Since $C$ corresponds to a local maximum, $C$ has shortening moves to both sides.    By the re-ordering principle, we can assume that the triangles corresponding to these shortening moves are $T_{j}$ and $T_{j+1}$.  

If $T_j$ and $T_{j+1}$ are disjoint or share a vertex, one can check that the new ordering obtained by interchanging $T_j$ and $T_{j+1}$,  
$$
O'(T):  T_1,...,T_{j-1},T_{j+1},T_j,T_{j+2},....,T_n
$$  
is also good, but $w_{O'}(T)<w_O(T)$.

\begin{figure}[h]
    \centering
    \includegraphics[width=0.8\textwidth]{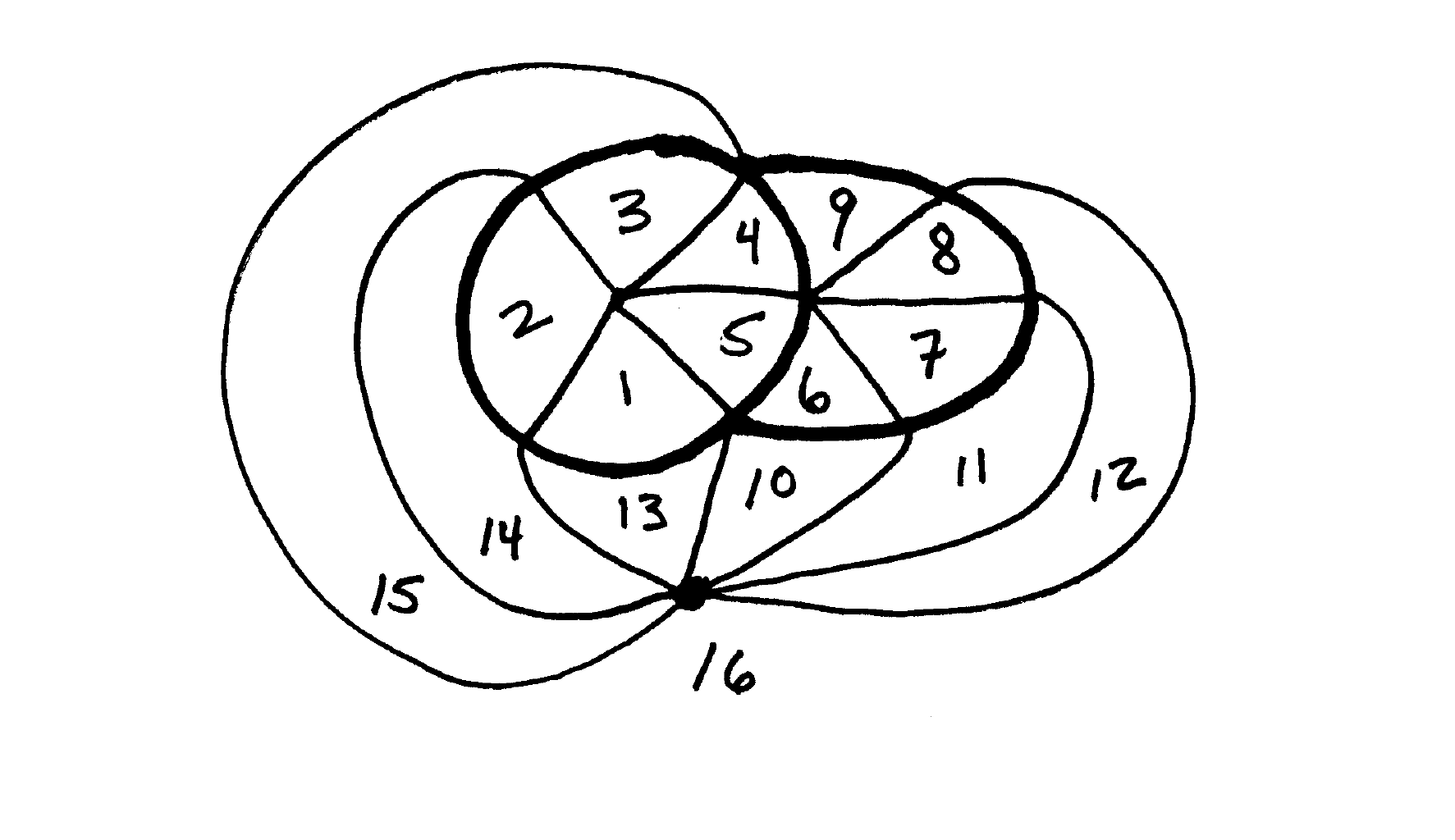}
    \caption{ $\partial I_5$ and $\partial I_9$ are adjacent local minima}
    \label{fig4}
\end{figure}

This contradicts the hypothesis that $O(T)$ is thin, hence $C$ is an unstable geodesic, proving part 1 of the theorem.    To conclude the proof, we consider what happens between two local minima:

\begin{Def} Suppose $T$ with order $O(T)$ is in thin position.   Suppose $\partial I_i$ and $\partial I_k$ are cycles in in $T$ corresponding to local minima of $O(T)$ such that the ordered list \\
$(|\partial I_i|, |\partial I_{i+1}|,....,|\partial I_{k}|)$ has a single local maximum.   Then $\partial I_i$ and $\partial I_k$ are {\it adjacent} minima in $T$.    If  the ordered list $(|\partial I_1|, |\partial I_{2}|,....,|\partial I_{i}|)$ has a single local maximum not at $|\partial I_{i}|$ we say that  $\partial I_i$ is adjacent to the empty geodesic (see Figure \ref{fig4}).
\end{Def}

Now suppose to the contrary that some cycle corresponding to a minimum of $O(T)$ at $|\partial I_i|$ has a shortening move.   Then either there exists one such move corresponding to a triangle lying between adjacent minima, or $\partial I_i$ is adjacent to the empty geodesic and there is a shortening move for  $\partial I_i$ in the disk $I_i$.   

Suppose $|\partial I_i|$ and $|\partial I_k|$ are adjacent minima in $T$, corresponding to cycles $C_i$ and $C_k$.   Assume the single maxima between them occurs at $C_j$.    Suppose $C_i$ has a shortening move across $T_h$ and $T_h$ is contained in the region between $C_i$ and $C_k$.   By the  re-ordering principle, we can re-order the triangles between $C_i$ and $C_k$, without increasing the width, so that $h=i+1$.   Since the addition of each triangle exactly increases or exactly decreases the length of the disk boundary by $1$,  the number of triangles between $C_i$ and $C_k$ is exactly $(|C_j|-|C_i|)+(|C_j|-|C_k|) $.   When we re-order the triangles so that $h=i+1$, the maximum length achieved is at least one smaller than $|C_j|$, hence the overall width is smaller than $O(T)$.  This contradicts thinness of $O$, hence $C_i$ cannot  have a shortening move across $T_h$ with $T_h$  contained in the region between $C_i$ and $C_k$.   

Suppose $\partial I_i$ is adjacent to the empty geodesics and there is a shortening move for  $\partial I_i$ in the disk $I_i$.   By the reordering principle we can reorder the triangulation so that the triangle associated to the shortening move is $T_i$.    However the width of this reordering is lower than that of the original ordering, contradicting thinness.      

Hence no cycle corresponding to a minimum has any shortening move, hence all such cycles are stable geodesics, as required.

\section {Width and Hamiltonian cycles}

A theorem of H. Whitney gives sufficient conditions for $T$ to contain a Hamiltonian cycle.  We examine the relation between the existence of a Hamiltonian cycle and bridge position for $T$.      Recall that $T$ with order $O(T)$ is in bridge position if the ordered list $(|\partial I_1|, |\partial I_2|,....,|\partial I_{n-1}|)$ has a single local maximum and no local minima.

\begin{Thm}{\cite{W}} If every cycle of length three in $T$ is the boundary of a triangle in $T$, then $T$ has a Hamiltonian cycle.
\end{Thm}

\begin{Thm}$T$ has a Hamiltonian cycle if and only if $T$ has an order $O(T)$ so that $T$ with order $O(T)$ is in bridge position.

\end{Thm}
\medskip

{\bf Proof}

Suppose $T$ has a Hamiltonian cycle. Let $D_1$ and $D_2$ be the two disks (thought of in $S^2$) defined by the Hamiltonian cycle.   Let $\Gamma_i$ be the graph dual to $T$ in $D_i$.   Since all the vertices of $T$ lie on $\partial D_i$, $\Gamma_i$ is a tree.     Construct the desired (good) order $O(T)$ by constructing $\Gamma_1$ from a root to the leaves, and then reversing the process for $\Gamma_2$.

\medskip

Conversely suppose the ordering $O(T)$ on $T$ has a unique local maximum and no local minima.  Suppose $|\partial I_j|$ is the unique local maximum for $O(T)$. Then $\partial I_j$ is a Hamiltonian cycle for $T$.  

\begin{Coro} If every cycle of length three in $T$ is the boundary of a triangle in $T$ then $T$ has an  order $O(T)$ so that $T$ with order $O(T)$ is in {\it bridge position}.
\end{Coro}

\section {When thin equals bridge}

\begin{Thm}

Let $T$ be a triangulation of the $2$-sphere.    Suppose $T$ with order $O(T)$ is in both thin position and bridge position.   Then $T$ is the tetrahedral triangulation of $S^2$.  
\end{Thm}

\medskip

{\bf Proof}

Suppose $T$ with order $O(T)$ is in both thin position and bridge position. Let $I_k$ be the disk such that $|\partial I_k|$ realizes the single local maximum of $O(T)$.   Let $J_k=S^2-I_k$.

By Theorem 11, $\partial I_k$ is a Hamiltonian cycle in the $1$-skeleton of $T$.  By Theorem 6,  $\partial I_k$ is an unstable geodesic, so cannot have disjoint, or 1-point intersecting, shortening moves in $I_k$ and $J_k$.  

\begin{figure}[h]
    \centering
    \includegraphics[width=0.8\textwidth]{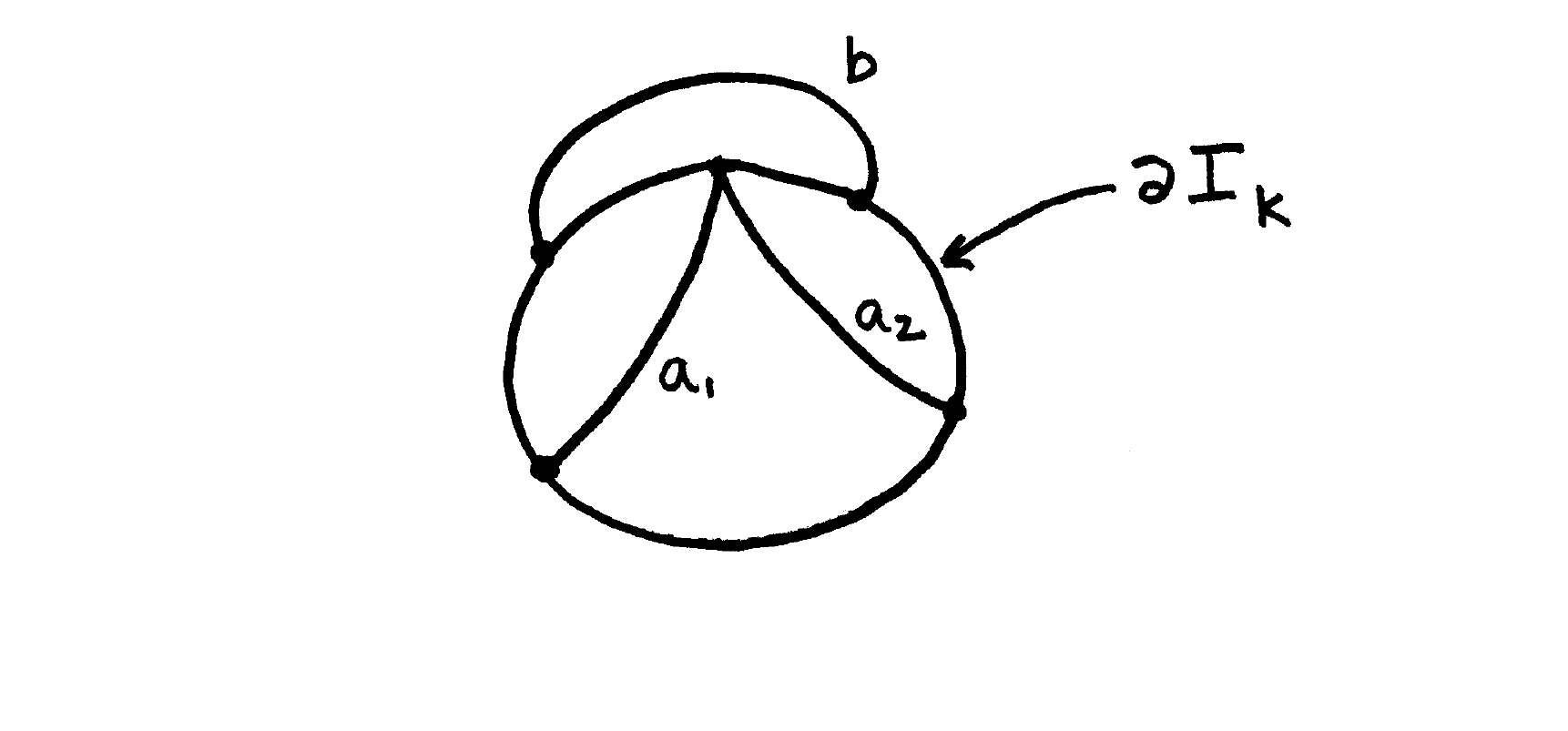}
    \caption{outermost arcs}
    \label{fig5}
\end{figure}

Each of $I_k$ and $J_k$ have at least two distinct outermost arcs, else a single arc which is outermost to both sides.  Each outermost arc corresponds to a shortening move, since there are no vertices in the interior of $I_k$ or $J_k$.    Let $a_1$ be an outermost arc of $I_k$, $b$ be an outermost arc of $J_k$.    Then the endpoints of $a_1$ and $b$ must be nested on $\partial(I_k)=\partial(J_k)$, else there will be shortening moves corresponding to $a$ and $b$ which are disjoint (see Figure \ref{fig5}).    

This nesting must hold for all possible pairs of outermost arcs in $I_k$ and $J_k$.   Suppose $a_2$ is a distinct outermost arc of $I_k$.    Since it must also have nested endpoints with $b$, it must share exactly one endpoint with $a_1$.  An additional outermost arc of $J_k$ will have to nest both with $a_1$ and with $a_2$, forcing it to coincide with $b$.  Hence $J_k$ has exactly one outermost arc, hence $I_k$ has exactly one outermost arc, and the theorem follows.       

Note that the tetrahedral triangulation has three length four unstable geodesics, similar to the smooth case of an ellipsoid with three distinct radii.

\section {When thin does not equal bridge}

\begin{Thm}

Let $T$ be a triangulation of the $2$-sphere.    Suppose $T$ with order $O(T)$ is in thin position but not bridge position.   Then $T$ has at least three distinct geodesics.

\end{Thm}

\medskip

{\bf Proof}

Suppose $T$ with order $O(T)$ is in thin position but not bridge position. Then the ordered list $(|\partial I_1|, |\partial I_2|,....,|\partial I_{n-1}|)$ has at least two local maxima, say at  $|\partial I_i|$ and $|\partial I_k|$,  and at least one local minima, say at $|\partial I_j|$.  Hence $T$ has at least two unstable geodesics, $\partial I_i$ and $\partial I_k$ and one stable geodesic, $\partial I_k$.   While distinct, they may overlap in paths.      

\begin{Coro} 
Let $T$ be a triangulation of the $2$-sphere.  Then $T$ has at least three distinct geodesics. 
\end{Coro}

\medskip

{\bf Proof}

$T$ is either the tetrahedral triangulation or there exists $O(T)$ such that $T$ with order $O(T)$ is in thin position but not bridge position.  The result follows from our observation on the tetrahedral triangulation and from the previous theorem.

\bigskip
By carefully considering regions between stable geodesics, we can improve this result, to obtain three distinct stable geodesics except in two cases.    We accomplish the needed details for this in the next section.

\section {Three geodesics revisited}

We begin with a theorem giving a precise description of the region between adjacent minima in a triangulation in thin position.   The result yielding three stable geodesics (in most cases) appears as a corollary.    

\begin{figure}[h]
    \centering
    \includegraphics[width=0.8\textwidth]{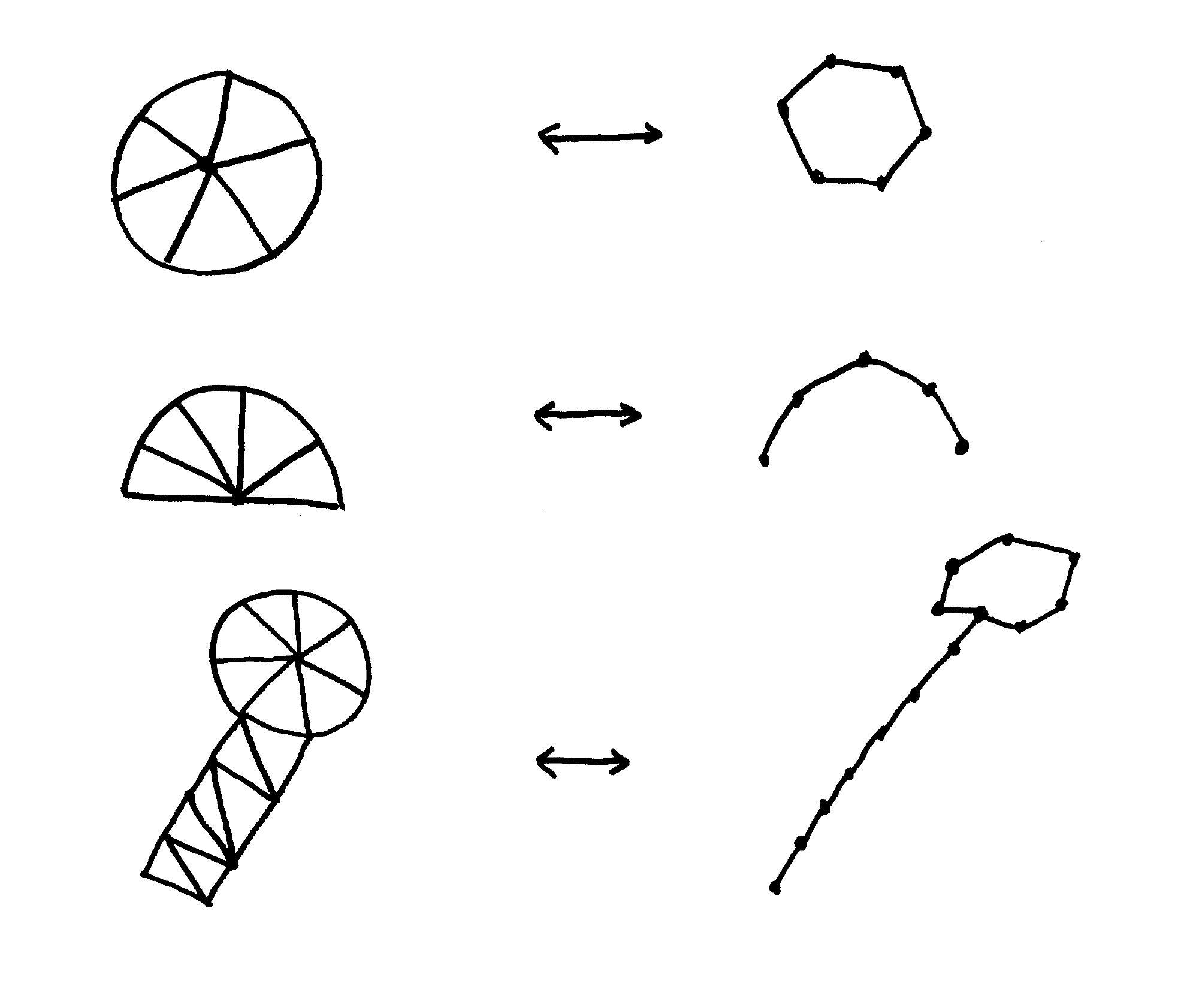}
    \caption{wheel, fan, planar lollipop, and their dual graphs}
    \label{fig6}
\end{figure}

\begin{Def} Let  $D$ be a triangulated disk, and let $\Gamma_T$ denote the dual graph to the triangulation.  A triangulated disk $D$ is a {\it wheel}  if $\Gamma_T$ is a cycle.    A triangulated disk $D$ is a {\it planar lollipop} if $\Gamma_T$ is isotopic to a cycle with an antenna attached.  A triangulated disk $D$ is a {\it fan} if $\Gamma_T$ is isotopic to an arc, and the two triangles corresponding to the endpoints of the arc share a vertex, the {\it distinguished vertex},  in $T$ (see Figure \ref{fig6}).  
\end{Def}

\begin{Thm}

Let $T$ be a triangulation of the $2$-sphere.   Assume $T$ with order $O(T)$ is in thin position.  Assume every cycle of length three bounds a triangle.  Suppose $\partial I_i$ and $\partial I_k$ are stable geodesics corresponding to adjacent minima in $T$.  Then $\partial I_i$ and $\partial I_k$ (or $\partial I_i$ alone) define a subdisk $G$ of the 2-sphere, with induced triangulation which is a wheel, a fan or a planar lollipop.   If $\partial I_i$ is adjacent to the empty geodesic, $G$ is a wheel.   

\end{Thm}

{\bf Proof} 

We begin with the following Claim:

\begin{Claim} 
Let $D$ be a triangulated disk with a (good) ordering $O(T)$ $T_1, T_2,...,T_m$  such that the ordered list $(|\partial I_1|, |\partial I_2|,....,|\partial I_{m}|)$ has a single local maximum at $|\partial I_{k}|, k\neq m$.   If $O(T)$ is thin, then $D$ is a wheel.    
\end{Claim}

{\bf Proof}

Let  $\alpha=\partial I_k$; we know that $\alpha$ is an unstable geodesic.    Let $\Gamma'_T$ be the dual graph of $I_k$.  $\Gamma'_T$ is a tree by Claim 7.  The disk $I_{k+1}$ is obtained from $I_k$ by adding the single disk $I_{k+1}$.   The effect of this addition on $\Gamma'_T$ is to add a single, 2-valent vertex, changing $\Gamma'_T$ from a tree to a graph $\Gamma''_T$ with a single cycle.    If  $\Gamma''_T$ has any 1-valent vertices, these correspond to shortening moves for $\alpha$ which are disjoint from $T_{k+1}$.   Hence the tree $\Gamma'_T$ can have at most two leaves, hence it must have exactly two leaves, both of which are connected to the new vertex corresponding to $T_{k+1}$ in $\Gamma''_T$.   Hence $I_{k+1}$ is a wheel.    If $I_{k+1}$ is a wheel with exactly three spokes, then no additional shortening move of the boundary is possible without violating the simplicial structure, so $m=k+1$, and $D$ is a wheel, as required.   Suppose $I_{k+1}$ is a wheel with strictly more than three spokes, and suppose $(k+1)<m$.  Let $\Gamma'''_T$ be the dual graph of $I_{k+2}$.  $\Gamma'''_T$ retracts onto a theta curve, with one loop of the theta curve a cycle of length at least $4$, corresponding to the wheel $I_{k+1}$, and one loop a cycle of length $3$.   We can re-order the triangles in $I_{k+2}$ to complete the length 3 cycle first, reducing the width of the triangulation in $I_{k+2}$ and hence in $D$, contradicting thinness of $O$.   Hence if $I_{k+1}$ is a wheel with strictly more than three spokes then $(k+1)=m$ and $D=I_{k+1}$, which is a wheel as required.    
\bigskip

\begin{figure}[h]
    \centering
    \includegraphics[width=0.8\textwidth]{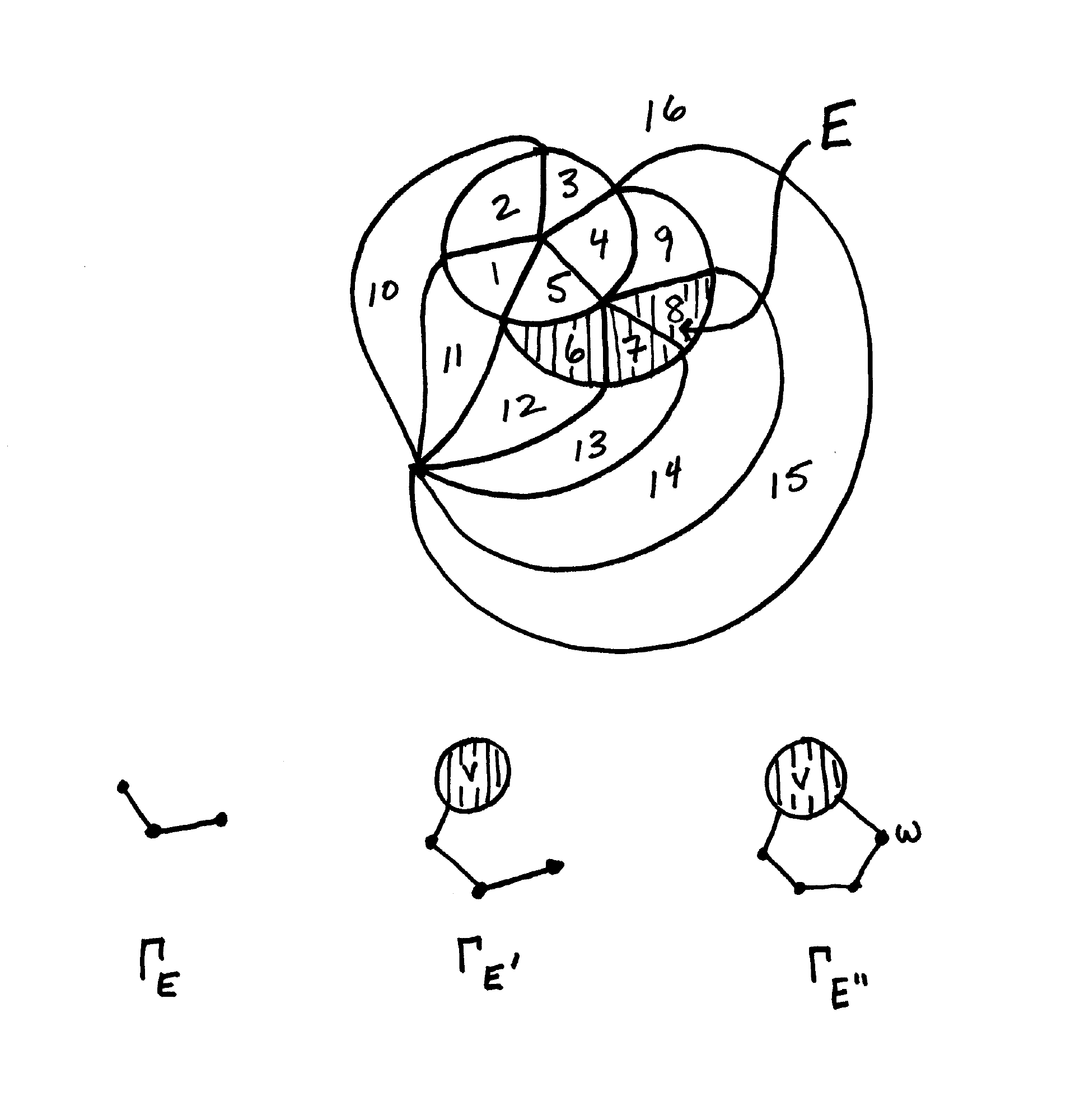}
    \caption{$E=(T_{6}\cup{T_{7}}\cup{T_8})$}
    \label{fig7}
\end{figure}

We continue with the proof of the Theorem:

Now assume $\partial I_i$ and $\partial I_k$ are adjacent stable geodesics.  Let $\alpha=\partial I_j$ be the unstable geodesic corresponding to the maximum that lies between $|\partial I_i|$ and $|\partial I_k|$.   Let $E$ consist of the triangles in $T$  between $|\partial I_i|$ and $|\partial I_j|$  , i.e., $E=(T_{i+1}\cup{T_{i+2}}\cup...\cup{T_j})$.   Let $\Gamma_E$ be the dual of $E$.  Let ${\Gamma_E}'$ be constructed from $\Gamma_E$ by adding a vertex $v$ corresponding to the disk $I_i$ and an edge for each triangle in $E$ sharing an edge with $\partial{I_i}$ (see Figure \ref{fig7}).    

Note that since each triangle in $\{T_{i+1},{T_{i+2}},...,{T_j}\}$ increases the length of the boundary of $I_i$ while leaving the homeomorphism type unchanged, ${\Gamma_E}'$ is a tree.     

Adding the triangle $T_{j+1}$ to $I_i\cup(T_{i+1}\cup{T_{i+2}}\cup...\cup{T_j})$  corresponds to adding a single bi-valent vertex $w$ to ${\Gamma_E}'$; call this new induced dual graph ${\Gamma_E}''$ .  

Recall $\alpha$ is an unstable geodesic.  Adding $w$ to ${\Gamma_E}'$ is a shortening move on $\alpha$.  A leaf of ${\Gamma_E}'$ which corresponds to a triangle in $E$ also corresponds to a shortening move on $\alpha$, hence the addition of  $T_{j+1}$  must eliminate all leaves of ${\Gamma_E}'$ which correspond to triangles in $E$, else there will be disjoint shortening moves on opposite sides of $\alpha$, a contradiction.    Hence ${\Gamma_E}'$ is a tree with at most two leaves corresponding to triangles in $E$, hence ${\Gamma_E}''$ is the dual of a wheel, a fan or a planar lollipop, with one additional vertex $v$ appended.    

Our goal is to show that $j=i+1$ or $j=k-1$ (or possibly both); that is, we would like to see that either we cannot add any (boundary reducing) further triangles to $I_{j+1}$ without violating thinness, and so we are done, or else we arrived at $\alpha$ after adding only a single triangle to $I_i$.   In the second case, we achieve the desired result by working backwards from the disk $S^{2}-(I_k)$.    

So assume ${\Gamma_E}''$ is the dual of a wheel, a fan or a planar lollipop, and assume that $j>i+1$.
Suppose also that $j<k-1$.   Then there is at least one additional triangle $U$ not contained in  $I_{k+1}$ which lies in $E$, and adding that triangle to   $I_{k+1}$ must decrease the length of the boundary.     We proceed by inspection to show that in every case, we can reduce the width of $O(T)$, violating thinness.    We examine the case when ${\Gamma_E}''$ corresponds to a planar lollipop; the others are similar.   

\begin{figure}[h]
    \centering
    \includegraphics[width=0.8\textwidth]{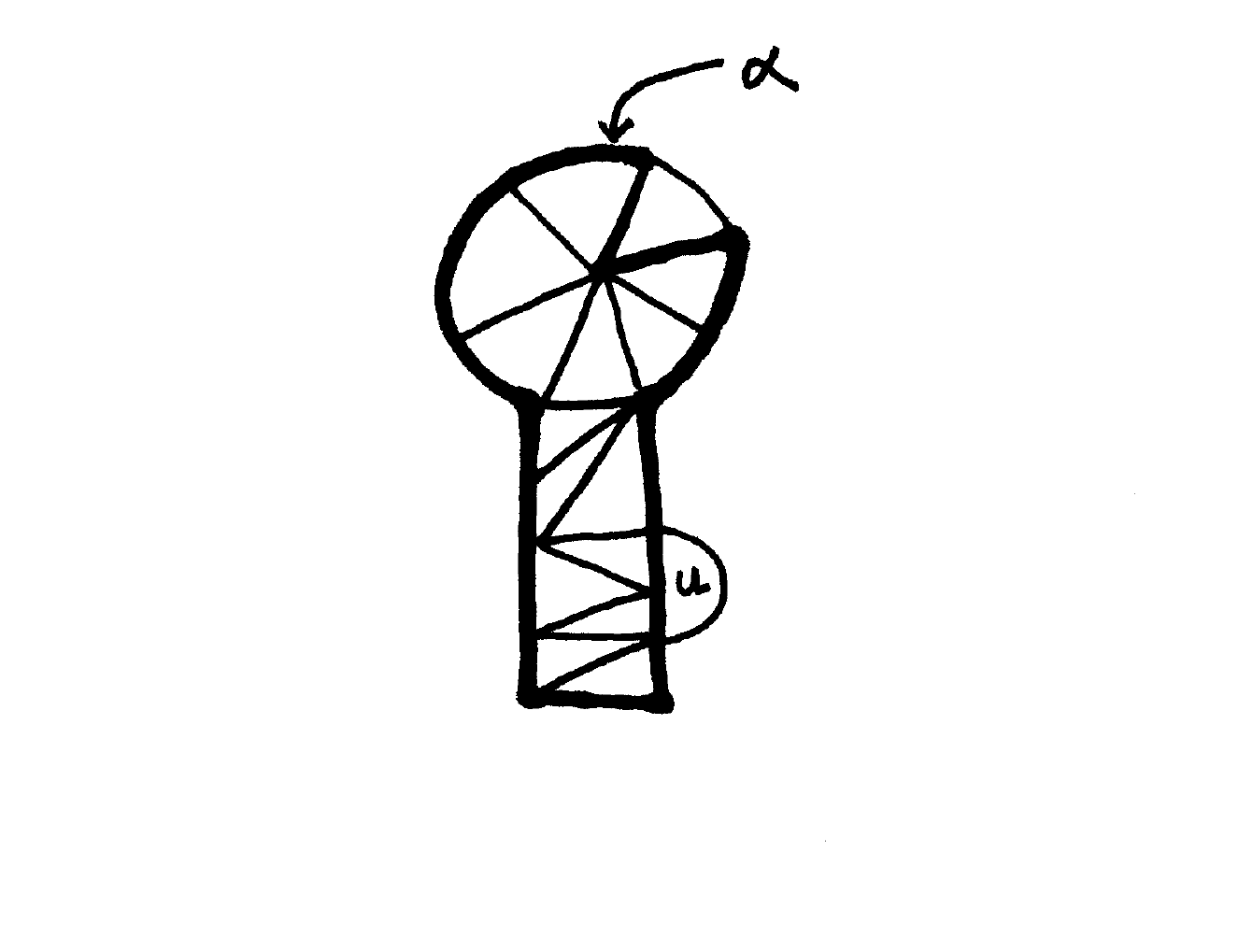}
    \caption{$U$ is adjacent only to the ``stick''}
    \label{fig8}
\end{figure}

 Note that since every cycle of length three bounds a triangle by assumption, $U$ cannot be adjacent only to the ``pop" section of the lollipop, so we need only consider the possibilities that is it adjacent to the ``pop'' and the stick, the stick alone, or the stick and the boundary of $I_i$.    In each case we observe that the addition of $U$ creates two disjoint shortening moves for $\alpha$ on opposite sides, a contradication to the assumption that $O(T)$ is thin (see one case in Figure \ref{fig8}).    
 
Since the existence of $U$ contradicts the hypothesis that $O(T)$ is thin, we conclude that $U$ cannot exist, hence $j=k-1$ and $E$ has the desired form.

\medskip

\begin{figure}[h]
    \centering
    \includegraphics[width=0.8\textwidth]{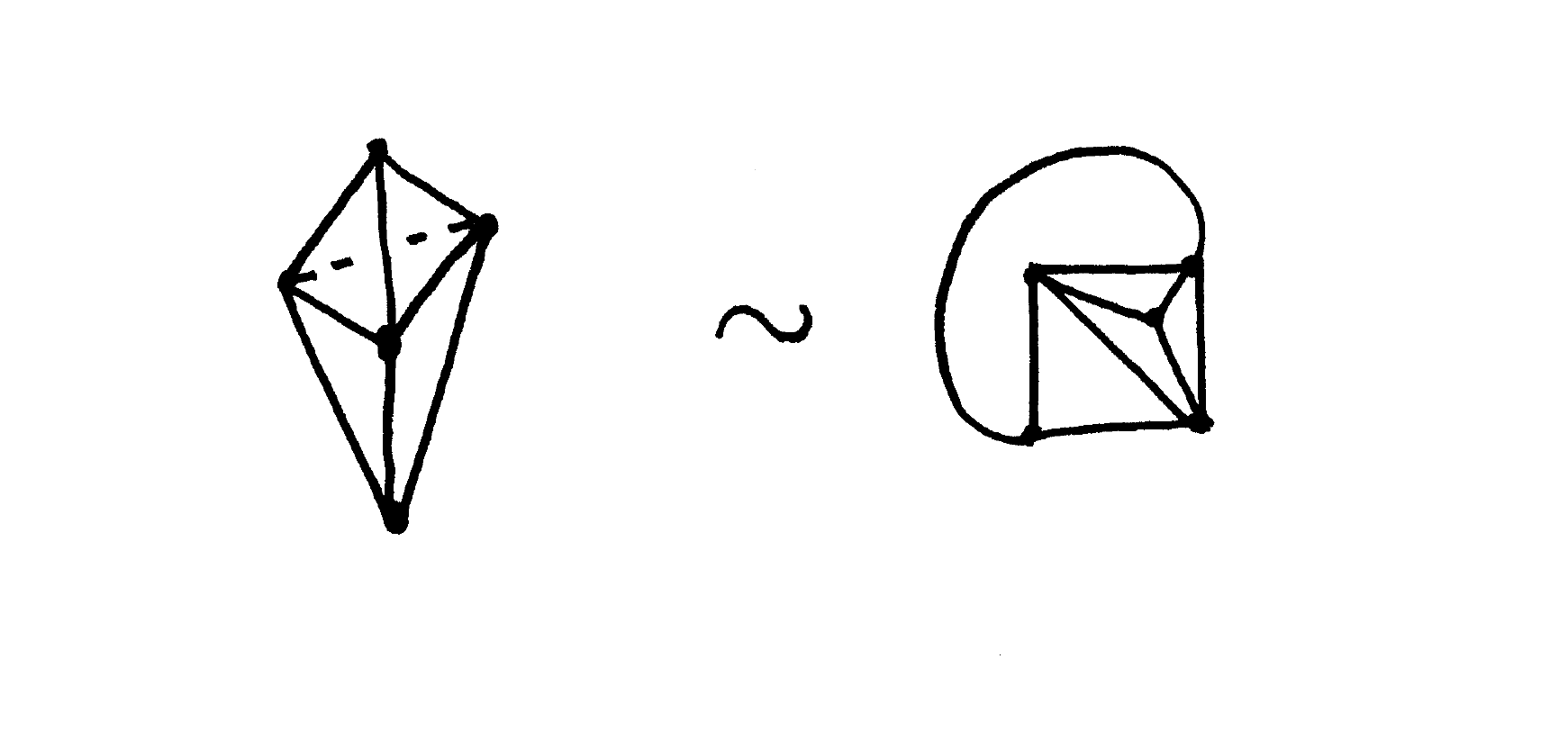}
    \caption{double tetrahedral triangulation}
    \label{fig9}
\end{figure}

\begin{Coro} 

Let $T$ be a triangulation of the $2$-sphere.   Then either $T$ is the tetrahedral triangulation, or the `` double tetrahedral'' triangulation (see Figure \ref{fig9}) obtained by attaching two tetrahedra along a single face,   or $T$ has at least three distinct stable  geodesics. 
\end{Coro}

{\bf Proof}  

Assume $T$ is not the tetrahedral triangulation.  Suppose $T$ with order $O(T)$ is in thin position.  Then $O(T)$ is not bridge position.    If $T$ with $O(T)$ has at least three local minima, by Theorem 6 they each correspond to a (distinct) stable geodesic and we are done.       Hence we need to consider the two cases $O(T)$ has exactly one local minimum and $O(T)$ has exactly two local minima.

Case 1:  Assume  $O(T)$ has exactly one local minimum.   Then, by Claim 18, the unique local minimum  for $O(T)$ splits the sphere into two wheels, $W$ and $V$.   $W$ and $V$ have the same number of spokes.  If the number of spokes in each wheel is three, then the triangulation is the double tetrahedral triangulation and we are done.
 
\begin{figure}[h]
    \centering
    \includegraphics[width=0.8\textwidth]{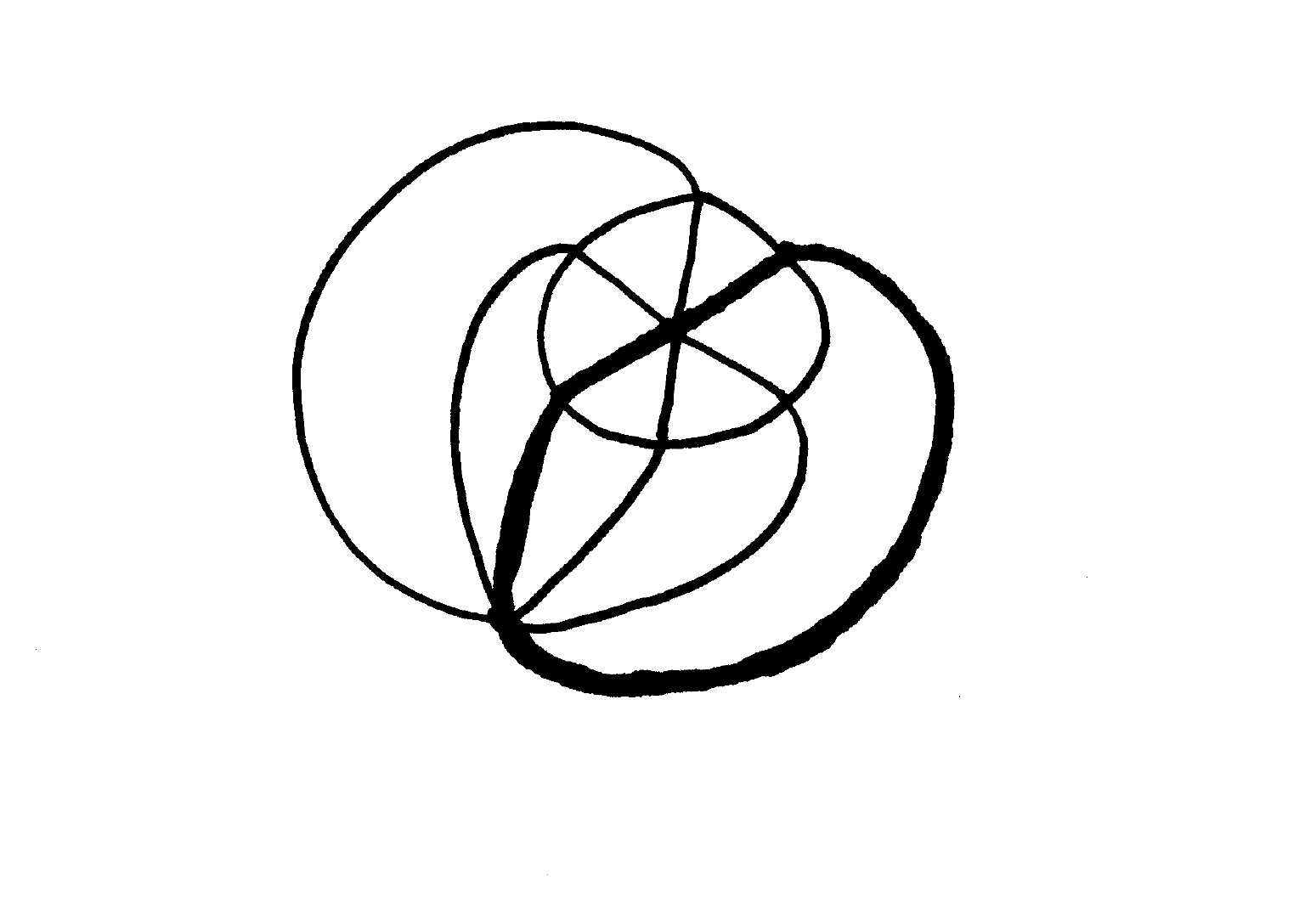}
    \caption{a stable geodesic through both hubs}
    \label{fig10}
\end{figure}

Suppose the number of spokes is at least four.   Then we can find (at least) two additional stable geodesics by constructing length four paths that contain the hubs of $V$ and $W$, including non-adjacent spokes in each wheel.   As there will be at least two such paths (see Figure \ref{fig10}), and these paths are stable geodesics, the result follows.

Case 2:  Assume $O(T)$ has exactly two local minima.   Then the two local minima correspond to distinct stable geodesics $\alpha$ and $\beta$; we need to find a third.    

If every length three geodesic in $T$ bounds a triangle, then by theorem 17, the region between the local minima, $E$, is a disk, and the triangulation restricted to $E$ is a wheel, a lollipop or a fan.   The triangulation in the complement of $E$ is two wheels partly attached along their rims.    If $E$ is a wheel or a lollipop, it contains a vertex of $T$ in its interior, and the link of that vertex is a stable geodesic distinct from $\alpha$ and $\beta$.   Suppose $E$ is a fan.   Let $v$ be the distinguished vertex of $E$.  Then $E$ is attached to one of the complementary wheels along two edges incident to $v$.    The link of $v$ is again a stable geodesic,  forming the boundary curve of a wheel with hub at $v$.     

Assume there exists a length three geodesic $\gamma$ in $T$ which does not bound a triangle.   Theorem 17 works as before unless $\gamma$ lies in the disk $E$ between $\alpha$ and $\beta$. In that case $\gamma$ is distinct from $\alpha$ and $\beta$, and provides the third stable geodesic we are seeking.

\begin{flushright}
Abigail Thompson\\ Department of Mathematics\\
University of California\\ Davis,
CA 95616\\ e-mail: thompson@math.ucdavis.edu\\
\end{flushright}


\begin{thebibliography}{HHH}

\bibitem{B} Bing, R. H. {\it Some aspects of the topology of 3-manifolds related to the
Poincare conjecture},  Lectures on Modern Mathematics, Volume II , Chapter 3, 93-128. John Wiley and Sons, 1964.


\bibitem{G} Grayson, Matthew A. {\it
Shortening embedded curves} 
Ann. of Math. (2) 129 (1989), no. 1, 71�111. 

\bibitem{H} Hempel, John {\it 3-Manifolds}, Annals of Mathematics Studies, vol. 86, Princeton University Press, 1976.

\bibitem{K} Klingenberg, Wilhelm
{\it Lectures on closed geodesics}
Third edition. Mathematisches Institut der Universit�t Bonn, Bonn, 1977. 210 pp.

\bibitem{L-S} L. Lusternick and L. Schnirelmann {\it Sur le probl\'{e}me de trois g\'{e}od\'{e}siques ferm\'{e}es sur les surfaces de genre $0$}, C.R. Acad. Sci. Paris 189 (1929), 269-271.

\bibitem{W} Whitney, H. {\it A theorem on graphs}, Ann. of Math. (2) 32 (1931), no. 2, 378-390.

%
 \end{thebibliography}
\end{document}